\documentclass[12pt]{article}
\usepackage{latexsym}
\usepackage{amsfonts,amsmath,amssymb}
\usepackage{epsfig}
\usepackage{francais}

\newtheorem{defi} {D{\'e}finition}
\newtheorem{theo} {Th{\'e}or{\`e}me}   
\newtheorem{prop} {Proposition}

\newtheorem{lem} {Lemme}

\begin{document}
    
     \bibliographystyle{abbrv}

  \begin{center}  
    
    \Huge  \bf

Cohomologie {\'e}quivariante des vari{\'e}t{\'e}s de Bott-Samelson.

\end{center}

\rm
\normalsize

 \vspace {0,25 cm}

 \begin{center}

\large

Matthieu WILLEMS

\end{center}

\rm
\normalsize

\vspace {1 cm}

Le but de ce texte est de donner explicitement les valeurs des
restrictions aux points fixes d'une base de la cohomologie
{\'e}quivariante des vari{\'e}t{\'e}s de Bott-Samelson. Ce calcul nous permet de
retrouver la structure de la cohomologie ordinaire de ces vari{\'e}tes. Il
permet {\'e}galement de retrouver les formules de S. Billey \cite{kl} pour la
restriction aux points fixes d'une base de la cohomologie {\'e}quivariante
des vari{\'e}t{\'e}s de Schubert. 

Je remercie Mich{\`e}le Vergne de m'avoir conseill{\'e} de regarder ce
probl{\`e}me et Alberto Arabia de m'avoir fait comprendre la structure de
ces vari{\'e}t{\'e}s et de leurs cohomologies.

\section{Cohomologie des vari{\'e}t{\'e}s de Bott-Samelson}

Soient $G$ un groupe de Lie semi-simple complexe et connexe,
$\mathfrak{g}$ son alg{\`e}bre de Lie, $\mathfrak{h}$ une sous-alg{\`e}bre de Cartan de $\mathfrak{g}$,
$\mathfrak{b}$ une sous-
alg{\`e}bre de Borel de $\mathfrak{g}$ contenant $\mathfrak{h}$, et $H$ et $B$ les sous-groupes de $G$ 
correspondant. Soient $T \subset K$ les formes r{\'e}elles de $H$ et $G$ 
respectivement. On pose $X=G/B=K/T$. On notera $e$ l'{\'e}l{\'e}ment neutre
 de $K$. Soit $W$ le groupe de Weyl 
correspondant {\`a} $(\mathfrak{g},\mathfrak{h})$, engendr{\'e} par les r{\'e}flexions simples $\{ r_{i} \}_{1 \leq i \leq r}$. On note $\{\alpha_{i} \}_{1 \leq i \leq r}
 \in \mathfrak{h}^*$ les racines simples correspondantes. On notera
$\Delta_{+} \subset\mathfrak{h}^*$ les racines positives et $\Delta_{-}$ les
racines n{\'e}gatives. On pose $S=S(\mathfrak{h}^*)$. On notera $v \leq w$
l'ordre de Bruhat dans $W$.

Soit $N$ un entier strictement positif. Consid{\'e}rons une suite de $N$ racines
simples $\mu_{1}$, ... , $\mu_{N}$ non n{\'e}cessairement distinctes. Pour $1 \leq i \leq N$, on note
$G_{i}$ le sous-groupe ferm{\'e} connexe de $G$ d'alg{\`e}bre de Lie
$\mathfrak{g_{\mu}}_{i}\oplus \mathfrak{h}\oplus
\mathfrak{g_{-\mu}}_{i}$ et on pose  $K_{i}=G_{i}\bigcap K$ . On d{\'e}finit : 

$$\Gamma(\mu_{1}, ... ,\mu_{N})=K_{1} \times_{T} K_{2} \times_{T} ... \times_{T}
K_{N}/T,$$
comme l'espace des orbites de $K_{1} \times K_{2} \times ... \times K_{N}$
sous l'action de $T \times T  \times ... \times T$ d{\'e}finie par : 

$$(t_{1}, t_{2}, ... , t_{N})(k_{1}, k_{2}, ... , k_{N}) = 
(k_{1}t_{1},t_{1}^{-1} k_{2}t_{2}, ... ,t_{N-1}^{-1} k_{N}t_{N}).$$

On notera $[k_{1}, k_{2}, ... , k_{N}]$ la classe de $(k_{1}, k_{2},
... , k_{N})$ dans $\Gamma(\mu_{1}, ... ,\mu_{N})$. On notera
$k_{\mu_{i}}$ un repr{\'e}sentant quelconque de la reflexion $r_{i}$ de
$N_{K_{\mu_{i}}}(T)/T$. Dans la suite, on
notera $\Gamma(\mu_{1}, ... ,\mu_{N})$ par $\Gamma$. On munit $\Gamma$
de sa structure complexe canonique d{\'e}finie dans \cite{ki}.

On d{\'e}finit une action de $T$ sur $\Gamma$ par : 

$$t[\mu_{1}, ... ,\mu_{N}]=[t\mu_{1}, ... ,\mu_{N}].$$

On pose $\mathcal{E} =
\{0,1\}^N$. Pour $\epsilon \in \mathcal{E}$, on note $Y_{\epsilon}
\subset \Gamma$ l'ensemble des classes $[k_{1},
k_{2}, ... , k_{N}]$ qui v{\'e}rifient pour tout entier $i$ compris entre
$1$ et $N$ : 

$$\left\{ \begin{array}{ll} k_{i} \in T
 & si \hspace{0,15 cm} \epsilon_{i} =0, \\ 
 k_{i} \notin T
 & si\hspace{0,15 cm} \epsilon_{i} =1.
\end{array}\right.$$

On v{\'e}rifie imm{\'e}diatement que cette d{\'e}finition est bien compatible avec
l'action de $T^N$. On munit $\mathcal{E}$ d'une structure de groupe en
identifiant $\{0,1\}$ avec $\mathbb{Z}/2\mathbb{Z}$. Pour $\epsilon \in \mathcal{E}$, on note
$\pi_{+}(\epsilon)$ l'ensemble des entiers $i$ tels que
$\epsilon_{i}=1$ et $\pi_{-}(\epsilon)$ l'ensemble des entiers $i$
tels que $\epsilon_{i}=0$. On pose $l(\epsilon) =
{\rm card}(\pi_{+}(\epsilon))$. On note $(i) \in \mathcal{E}$
l'{\'e}l{\'e}ment de $\mathcal{E}$ d{\'e}fini par $(i)_{j}=\delta_{i,j}$. Pour $\epsilon \in \mathcal{E}$, on pose
$\displaystyle{v_{i}(\epsilon) =\prod_{k \leq i, k \in \pi_{+}(\epsilon)
  }r_{\mu_{k}}}$, $v(\epsilon)=v_{l(\epsilon)}(\epsilon)$  et
$\alpha_{i}(\epsilon)=v_{i-1}(\epsilon)\mu_{i}$. Pour $i \leq j$, on
d{\'e}finit {\'e}galement $\displaystyle{v_{i}^j(\epsilon) =\prod_{i\leq k \leq j, k \in \pi_{+}(\epsilon)
  }r_{\mu_{k}}}$. Par convention, on pose $v_{0}(\epsilon)=e$ et donc
$\alpha_{1}(\epsilon)=\mu_{1}$. De plus, si $j<i$, on pose
$v_{i}^j(\epsilon)=e$. On d{\'e}finit un ordre sur $\mathcal{E}$ par : 

$$\epsilon \leq \epsilon' <=> \pi_{+}(\epsilon) \subset
\pi_{+}(\epsilon').$$

On d{\'e}montre alors facilement la proposition suivante : 

\begin{prop}
    
(i) Pour tout $\epsilon$ de $\mathcal{E}$, $Y_{\epsilon}$ est un
espace affine de dimension r{\'e}elle $2l(\epsilon)$. 

(ii) Pour tout $\epsilon$ de $\mathcal{E}$, $\overline{Y_{\epsilon}} = 
\coprod_{\epsilon' \leq \epsilon} Y_{\epsilon'}$

(iii) $\Gamma = \coprod_{\epsilon \in
  \mathcal{E}} Y_{\epsilon}$

(iv) Pour tout $\epsilon$ de $\mathcal{E}$, $Y_{\epsilon}$ est stable
par l'action de $T$.
   
\end{prop}

De plus, nous allons avoir besoin du lemme suivant :

\begin{lem}
    
(i) L'ensemble $\Gamma^T$ des points fixes de $\Gamma$ sous l'action de
$T$ est constitu{\'e} des $2^N$ points :

$$[k_{1}, k_{2}, ..., k_{N}], \hspace{0,15 cm} o\grave{u} \hspace{0,15 cm}
k_{i} \in \{ e, k_{\mu_{i}} \}.$$

On identifiera donc $\Gamma^T$ avec $\mathcal{E}$ en identifiant $e$
avec $0$ et $k_{\mu_{i}}$ avec $1$.

(ii) Soit $(\epsilon,\epsilon') \in \mathcal{E}^2$, alors : 

$$ \epsilon \in \overline{Y_{\epsilon'}} <=> \epsilon \leq \epsilon',$$

et dans ce cas : 

$$Pf(\epsilon, \overline{Y_{\epsilon'}})=(-1)^{l(\epsilon)}
\prod_{i \in \pi_{+}(\epsilon')}\alpha_{i}(\epsilon),$$
o{\`u} $Pf(\epsilon, \overline{Y_{\epsilon'}})$ d{\'e}signe le d{\'e}terminant de
l'action de $T$ sur l'espace tangent {\`a} $\overline{Y_{\epsilon'}}$ en
$\epsilon$ orient{\'e} par sa structure complexe d{\'e}duite de celle de
$\Gamma$.

\end{lem}

La d{\'e}composition $\Gamma=\coprod_{\epsilon \in
  \mathcal{E}}Y_{\epsilon}$ munit donc $\Gamma$ d'une structure de
  $CW$-complexe $T$ {\'e}quivariant o{\`u} toutes les cellules sont de
  dimension
 paire; de plus l'ensemble des points fixes de l'action de $T$ sur
  $\Gamma$ est discret. La proposition suivante est alors
  prouv{\'e}e dans \cite{kg} :

\begin{prop}

(i) La cohomologie $T$-{\'e}quivariante de $\Gamma^T$ s'identifie {\`a}
l'alg{\`e}bre des fonctions de $\mathcal{E}$ dans $S$, qu'on notera $F(\mathcal{E};S)$.

(ii) La restriction aux points fixes $i_{T}^*$ : $H_{T}^*(\Gamma)
 \rightarrow F(\mathcal{E};S)$ est injective.
    
(iii) La cohomologie $T$-{\'e}quivariante de $\Gamma$ est un $S$-module libre
qui admet comme base la famille $\{\hat{\sigma}_{\epsilon}\}_{\epsilon \in
\mathcal{E}}$ caract{\'e}ris{\'e}e par : 

$$\int_{\overline{Y_{\epsilon'}}}\hat{\sigma}_{\epsilon}=\delta_{\epsilon', \epsilon}.$$
   
\end{prop}

\begin{defi}
    
Pour $\epsilon \in \mathcal{E}$, on d{\'e}finit $\sigma_{\epsilon} \in
F(\mathcal{E};S)$
par :

$$\left\{ \begin{array}{ll} \sigma_{\epsilon}(\epsilon') = 
\displaystyle{\prod_{i \in \pi_{+}(\epsilon)}\alpha_{i}(\epsilon')}
 & si \hspace{0,15 cm} \epsilon \leq \epsilon' \\ 
 \sigma_{\epsilon}(\epsilon') = 0
 & sinon
\end{array}\right.$$

\end{defi}

On a alors la proposition suivante :

\begin{theo}
    
Pour tout $\epsilon \in \mathcal{E}$, on a : 

$$i_{T}^*(\hat{\sigma}_{\epsilon})=\sigma_{\epsilon}.$$
   
\end{theo}

\sc 
\begin{flushleft}
D{\'e}monstration :
\end{flushleft} 
\rm

En utilisant la formule de localisation et le lemme $1$, on
obtient pout tout  $\hat{\sigma} \in H_{T}^*(\Gamma)$ : 

$$\int_{\overline{Y_{\epsilon}}}\hat{\sigma}=(-1)^{l(\epsilon)}  \sum_{\epsilon' \leq
\epsilon}\frac{i_{T}^*(\hat{\sigma})(\epsilon')}{(-1)^{l(\epsilon')}\prod_{i
  \in \pi_{+}(\epsilon)}\alpha_{i}(\epsilon')}. \hspace{2 cm} (*)$$

Soit $\epsilon_{0} \in \mathcal{E}$, et soit
$\sigma'_{\epsilon_{0}}=i_{T}^*(\hat{\sigma}_{\epsilon_{0}})$. Montrons par
r{\'e}currence sur la longueur de $\epsilon$ que pour tout $\epsilon \in
\mathcal{E}$,
$\sigma'_{\epsilon_{0}}(\epsilon)=\sigma_{\epsilon_{0}}(\epsilon)$.
Grace {\`a} la
formule $(*)$ et {\`a} la caract{\'e}risation de $\hat{\sigma}_{\epsilon_{0}}$,
on d{\'e}montre facilement par r{\'e}currence sur $l(\epsilon)$ que si $\epsilon$ n'est pas plus grand que $\epsilon_{0}$, on
a bien $\sigma'_{\epsilon_{0}}(\epsilon)=0$. On peut donc se limiter
au cas o{\`u} $\epsilon_{0} \leq \epsilon$. Si $\epsilon=\epsilon_{0}$, la
formule $(*)$ et le fait que
$\int_{\overline{Y_{\epsilon_{0}}}}\hat{\sigma}_{\epsilon_{0}}=1$ nous
donne bien
$\sigma'_{\epsilon_{0}}(\epsilon_{0})=\sigma_{\epsilon_{0}}(\epsilon_{0})$.
  Soit $\epsilon > \epsilon_{0}$. On suppose
le r{\'e}sultat v{\'e}rifi{\'e} pour tout $\epsilon'$ de longueur strictement plus
petite que $\epsilon$, on applique la formule $(*)$ et le fait
que $\int_{\overline{Y_{\epsilon}}}\hat{\sigma}_{\epsilon_{0}}=0$  pour
obtenir : 

$$\sum_{\epsilon_{0} \leq  \epsilon' <
\epsilon}\frac{\prod_{i
  \in \pi_{+}(\epsilon_{0})}\alpha_{i}(\epsilon')}{(-1)^{l(\epsilon')}\prod_{i
  \in
  \pi_{+}(\epsilon)}\alpha_{i}(\epsilon')}+\frac{\sigma'_{\epsilon_{0}}(\epsilon)}{(-1)^{l(\epsilon)}\prod_{i
  \in \pi_{+}(\epsilon)}\alpha_{i}(\epsilon)}  = 0,$$

d'o{\`u} : 

$$\sum_{\epsilon_{0} \leq \epsilon' <
\epsilon}\frac{(-1)^{l(\epsilon')}}{\prod_{i
  \in
  \pi_{+}(\epsilon)\setminus \pi_{+}(\epsilon_{0})}\alpha_{i}(\epsilon')}+
\frac{(-1)^{l(\epsilon)}\sigma'_{\epsilon_{0}}(\epsilon)}{\prod_{i
  \in \pi_{+}(\epsilon)}\alpha_{i}(\epsilon)}  = 0.$$

Si on pose $\tilde{\epsilon}=\epsilon-(j)$, o{\`u} $j$ est le plus grand
{\'e}l{\'e}ment de $\pi_{+}(\epsilon) \setminus  \pi_{+}(\epsilon_{0})$, on a alors : 

$$\sum_{
 \tiny \begin{array}{ll}  \epsilon_{0} \leq \epsilon' <
\epsilon \\
\hspace{0,25 cm} \epsilon' \neq \tilde{\epsilon} 
\end{array}}
\frac{(-1)^{l(\epsilon')}}{\prod_{i
  \in
  \pi_{+}(\epsilon)\setminus \pi_{+}(\epsilon_{0})}\alpha_{i}(\epsilon')}-
\frac{(-1)^{l(\epsilon)}}{\prod_{i
  \in
  \pi_{+}(\epsilon)\setminus \pi_{+}(\epsilon_{0})}\alpha_{i}(\epsilon)}+
\frac{(-1)^{l(\epsilon)}\sigma'_{\epsilon_{0}}(\epsilon)}{\prod_{i
  \in \pi_{+}(\epsilon)}\alpha_{i}(\epsilon)}
  = 0.$$

En effet, comme $j$ est le plus grand {\'e}l{\'e}ment de $\pi_{+}(\epsilon) \setminus
\pi_{+}(\epsilon_{0})$, pour tout $i \in \pi_{+}(\epsilon) \setminus
\pi_{+}(\epsilon_{0})$,
$\alpha_{i}(\epsilon)=\alpha_{i}(\tilde{\epsilon})$ et de plus,
$l(\tilde{\epsilon})=l(\epsilon)-1$. Pour les m{\^e}mes raisons, on
s'aper{\c c}oit, en distinguant les termes qui ont un $1$ en $j${\`e}me position
et ceux qui ont un $0$ en $j${\`e}me position, que la premi{\`e}re somme est
nulle et on obtient alors bien $\sigma'_{\epsilon_{0}}(\epsilon)=
\prod_{i \in \pi_{+}(\epsilon_{0})}\alpha_{i}(\epsilon)$.

\vspace{0,5cm}

 On pose 
$\sigma_{i}=\sigma_{(i)}$ et 
$\hat{\sigma}_{i}=\hat{\sigma}_{(i)}$. On a alors : 

\begin{prop}
    
Pour tout $\epsilon \in \mathcal{E}$, on a : 

$$\hat{\sigma}_{\epsilon}=\prod_{i \in
  \pi_{+}(\epsilon)}\hat{\sigma}_{i}.$$

De plus, on a la formule de multiplication suivante : 

$$\left\{ \begin{array}{ll} \hat{\sigma}_{i}\hat{\sigma}_{\epsilon}=
\hat{\sigma}_{\epsilon+(i)}
 & si \hspace{0,15 cm} i\in \pi_{-}(\epsilon) \\ 
\displaystyle{ \hat{\sigma}_{i}\hat{\sigma}_{\epsilon}=\sigma_{i}(\epsilon)\hat{\sigma_{\epsilon}}
+\sum_{j<i,j\in
  \pi_{-}(\epsilon)}\frac{r_{j}\alpha_{j}^{i}(\epsilon)-\alpha_{j}^{i}(\epsilon)
 }{\mu_{j}}\hat{\sigma}_{\epsilon}}\hat{\sigma}_{j}
 & si \hspace{0,15 cm} i\in \pi_{+}(\epsilon)\end{array}\right.,$$

o{\`u} on a pos{\'e} $\alpha_{j}^{i}(\epsilon)=v_{j+1}^{i-1}(\epsilon)(\mu_{i})$.
   
\end{prop}

\sc 
\begin{flushleft}
D{\'e}monstration :
\end{flushleft} 
\rm

Par injectivit{\'e} de la restriction aux points fixes, il suffit de
d{\'e}montrer ces formules pour $\sigma_{i}$ et $\sigma_{\epsilon}$.

La premi{\`e}re formule se voit imm{\'e}diatement sur la d{\'e}finition des
$\sigma_{\epsilon}$.

Pour la deuxi{\`e}me formule, pour des raisons de degr{\'e} et d'ordre sur
$\mathcal{E}$, on sait que le produit $\sigma_{i}\sigma_{\epsilon}$ 
s'{\'e}crit sous la forme : 

$$\sigma_{i}\sigma_{\epsilon}=C_{i}\sigma_{\epsilon}+\sum_{j \neq i,j\in
  \pi_{-}(\epsilon)}C_{j}\sigma_{\epsilon+(j)}.$$

En {\'e}valuant en $\epsilon$, on trouve $C_{i}$ puis en {\'e}valuant en
$\epsilon+(j)$, on trouve : 
 $$C_{j}=(\prod_{k<j}r_{\mu_{k}})(\frac{r_{j}\alpha_{j}^{i}(\epsilon)-\alpha_{j}^{i}(\epsilon)
 }{\mu_{j}}).$$

 Or comme  $\frac{r_{j}\alpha_{j}^{i}(\epsilon)-\alpha_{j}^{i}(\epsilon)
 }{\mu_{j}}$ est un nombre, $C_{j}=(\prod_{k<j}r_{\mu_{k}})(\frac{r_{j}\alpha_{j}^{i}(\epsilon)-\alpha_{j}^{i}(\epsilon)
 }{\mu_{j}})=\frac{r_{j}\alpha_{j}^{i}(\epsilon)-\alpha_{j}^{i}(\epsilon)
 }{\mu_{j}}$ et on a bien la formule {\'e}nonc{\'e}e.

\vspace{0,5 cm}

En particulier, si pour $j<i$, on pose $a_{j,i}=n(\mu_{j}, \mu_{i})$, o{\`u}
$n(\alpha, \beta)$ d{\'e}signe le nombre de Cartan associ{\'e} aux racines
$\alpha$ et $\beta$, on a : 

\begin{prop}

$$  \hat{\sigma}_{i}^2=\alpha_{i}\hat{\sigma}_{i} -  
 \sum_{j<i}a_{j,i}\hat{\sigma}_{i}\hat{\sigma}_{j}.$$

\end{prop}

Grace {\`a} la d{\'e}composition cellulaire $\Gamma = \coprod_{\epsilon \in
  \mathcal{E}}Y_{\epsilon}$, une base de la cohomologie ordinaire de $\Gamma$ est
  aussi index{\'e}e par $\mathcal{E}$. On la notera $(x_{\epsilon})_{\epsilon \in
\mathcal{E}}$ et on notera $x_{i}=x_{(i)}$. Soit $v_{0}$ l'{\'e}valuation {\`a} l'origine
$H_{T}^*(\Gamma) \rightarrow H^*(\Gamma)$. Le r{\'e}sultat suivant est
  imm{\'e}diat :

\begin{lem}
    
L'{\'e}valuation {\`a} l'origine est l'homomorphisme d'anneaux d{\'e}fini par : 

$$v_{0}(f\hat{\sigma}_{\epsilon})=f(0)x_{\epsilon}.$$
   
\end{lem}

Grace {\`a} cet homomorphisme, on retrouve alors le r{\'e}sultat suivant
prouv{\'e} dans \cite{kf} : 

\begin{prop}

La cohomologie ordinaire de $\Gamma$ est engendr{\'e} par des {\'e}l{\'e}ments
$(x_{i})_{1 \leq i \leq N}$ de degr{\'e} $2$ soumis aux relations : 

$$x_{i}^2 + \sum_{j<i}a_{j,i}x_{i}x_{j}=0.$$

\end{prop}

\section{Applications aux vari{\'e}t{\'e}s de Schubert} 

Pour $w \in W$, on d{\'e}finit $X_{w}=\phi(BwB/B)$ o{\`u} $\phi$ est
l'isomorphisme canonique entre $G/B$ et $K/T$. La d{\'e}composition
$X=\coprod_{w\in W}X_{w}$ munit $X$ d'une structure de $CW$ complexe $T$ {\'e}quivariant o{\`u} toutes les cellules
sont de dimension paire et on a donc la
proposition suivante : 

\begin{prop}

(i) La cohomologie $T$-{\'e}quivariante de $X^T$ s'identifie {\`a}
l'alg{\`e}bre des fonctions de $W$ dans $S$, qu'on notera $F(W;S)$.

(ii) La restriction aux points fixes $i_{T}^*$ : $H_{T}^*(X)
 \rightarrow F(W;S)$ est injective.
    
(iii) La cohomologie $T$-{\'e}quivariante de $X$ est un $S$-module libre
qui admet comme base la famille $\{\hat{\xi}^{w}\}_{w \in W}$
caract{\'e}ris{\'e}e par : 

$$\int_{\overline{X_{w'}}}\hat{\xi}^{w}=\delta_{w',w}.$$
   
\end{prop}

On pose  $\xi^{w}=i_{T}^*(\hat{\xi}^{w})$. Soit $(w,v) \in W^2$. Soit 
$v=r_{i_{1}}\ldots r_{i_{l}}$ une d{\'e}composition r{\'e}duite de $v$. 
Pour $1 \leq j \leq l$, on d{\'e}finit un {\'e}l{\'e}ment $\beta_{j}$ de 
$S$ par $\beta_{j}=r_{i_{1}}\ldots 
r_{i_{j-1}}\alpha_{i_{j}}$. La formule suivante est prouv{\'e}e par
Sarah Billey dans \cite{kl} : 

\begin{prop}

    Soit $(w,v) \in W^2$, tels que $l(w)=m$ et $l(v)=l$. On a : 
   $$ \xi^{w}(v) = \sum\beta_{j_{1}}\ldots \beta_{j_{m}},$$
    o{\`u} la somme porte sur l'ensemble des entiers $1\leq j_{1} 
   < \ldots < j_{m} \leq l$ tels que $w=r_{i_{j_{1}}}\ldots 
   r_{i_{j_{m}}}$.
    
\end{prop}

Retrouvons cette formule gr{\^a}ce aux r{\'e}sultats pr{\'e}c{\'e}dents. Soit $w
\in W$ et soit  $w=r_{\mu_{1}} ...r_{\mu_{k}}$ une d{\'e}composition
r{\'e}duite de $w$. On d{\'e}finit une
application $g_{\mu_{1}, ..., \mu_{k}}$ de $\Gamma(\mu_{1}, ..., \mu_{k})$ dans $\overline{X_{w}}$ par
multiplication (i.e. $g([k_{1},... ,k_{k}])=k_{1}*... *k_{k}
\hspace{0.3cm}  [T]$). Le r{\'e}sultat suivant est prouv{\'e} dans  \cite{ki} :

\begin{prop}

L'application $g_{\mu_{1}, ..., \mu_{k}}$ est une d{\'e}singularisation de $\overline{X_{w}}$.
    
\end{prop}

Soit $w_{0}$ le plus grand {\'e}l{\'e}ment de $W$. Soit $w_{0}=r_{\mu_{1}} ...
r_{\mu_{N}}$ une d{\'e}composition r{\'e}duite de $w_{0}$. On pose $\Gamma=
\Gamma(\mu_{1}, ..., \mu_{N})$ et $g=g_{\mu_{1}, ..., \mu_{N}}$. La
proposition suivante va nous permettre de retrouver la proposition $7$
:

\begin{prop}

Soit $w \in W$, on a : 

$$g^*(\hat{\xi}^w)=\sum_{ \tiny \begin{array}{ll}\epsilon \in \mathcal{E}, l(\epsilon)=l(w) \\ 
\hspace{0,25 cm}  et \hspace{0,15 cm}  v(\epsilon)=w
\end{array}} \hat{\sigma}_{\epsilon},$$
    
\end{prop}

\sc 
\begin{flushleft}
D{\'e}monstration :
\end{flushleft} 
\rm

Pour d{\'e}montrer la proposition, il faut montrer : 

$$\int_{\overline{Y_{\epsilon}}}g^*(\hat{\xi}^w)=\delta_{v(\epsilon),w}
\delta_{l(\epsilon),l(w)}.$$

Pour des raisons de degr{\'e}, cette formule est {\'e}vidente quand
$l(\epsilon) > l(w)$. 

Quand $l(\epsilon)=l(w)$ et $v(\epsilon)=w$, comme
 $g_{/\overline{Y_{\epsilon}}} \hspace{0,3 cm} \overline{Y_{\epsilon}}  \rightarrow \overline{X_{w}}$
 s'identifie {\`a} l'application $g_{\mu_{i}, i \in
   \pi_{+}(\epsilon)}$, d'apr{\`e}s la proposition $8$, on a bien
 $\int_{\overline{Y_{\epsilon}}}g^*(\hat{\xi}^w)=\int_{\overline{X_{w}}}\hat{\xi}^w=1$.

Supposons d{\'e}sormais $l(\epsilon) \leq l(w)$ et $v(\epsilon) \neq
w$. Distinguons deux cas. 

Si $\epsilon$ correspond {\`a} une d{\'e}composition r{\'e}duite de
$v(\epsilon)$, par le m{\^e}me raisonnement que pr{\'e}c{\'e}demment, on a 
$\int_{\overline{Y_{\epsilon}}}g^*(\hat{\xi}^w)=\int_{\overline{X_{v(\epsilon)}}}\hat{\xi}^w=0$.

Si $\epsilon$ n'est pas une d{\'e}composition r{\'e}duite de $v(\epsilon)$, 
alors $g$ envoie $\overline{Y_{\epsilon}}$ dans
$\overline{X_{v(\epsilon)}}$ qui est de
dimension strictement plus petite que $\overline{Y_{\epsilon}}$ et donc
$\int_{\overline{Y_{\epsilon}}}g^*(\hat{\xi}^w)=0$.

\vspace {0,5 cm} 

De cette proposition, on d{\'e}duit que pour tous $w \in W$ et $\epsilon \in
 \mathcal{E}$, on a :

$$\xi^w(v(\epsilon))=\sum_{ \tiny \begin{array}{ll}\epsilon' \in \mathcal{E}, l(\epsilon')=l(w) \\ 
\hspace{0,25 cm}  et \hspace{0,15 cm}  v(\epsilon')=w
\end{array}}\sigma_{\epsilon'}(\epsilon),$$
ce qui nous redonne bien la proposition $7$, {\`a} l'aide du th{\'e}or{\`e}me $1$.

  \bibliography{biblio}
 \bibliographystyle{abbrv}

 \vspace{0,5 cm}

Universit{\'e} Paris 7, UFR de Math{\'e}matiques, case 7012

2, place Jussieu 75251 PARIS Cedex 05.

\vspace{0,2 cm}

\bf{e-mail} \rm : willems@math.jussieu.fr

\end{document}